\documentclass[a4paper,12pt]{amsart} 

\usepackage{amsmath,amssymb,amsxtra,amsrefs} 
\usepackage[all]{xy} 

\usepackage{hyperref} 
\usepackage[margin=3cm]{geometry} 
\newtheorem{theorem}{Theorem}

\theoremstyle{definition} 
\newtheorem{remark}[theorem]{Remark} 
\newtheorem{example}[theorem]{Example} 

\newcommand{\C}{\mathbb{C}} 
\newcommand{\R}{\mathbb{R}}
\newcommand{\Z}{\mathbb{Z}} 
\newcommand{\PP}{\mathbb{P}}

\newcommand{\chT}{\check{T}} 
\newcommand{\chK}{\check{K}} 

\newcommand{\hGamma}{\widehat{\Gamma}} 
\newcommand{\hI}{\widehat{I}} 
\newcommand{\hcI}{\widehat{\cI}}

\newcommand{\cY}{\mathcal{Y}} 
\newcommand{\cM}{\mathcal{M}} 
\newcommand{\cO}{\mathcal{O}} 
\newcommand{\cI}{\mathcal{I}} 

\newcommand{\bz}{\boldsymbol{0}} 

\newcommand{\frk}{\mathfrak{k}}

\newcommand{\iu}{\sqrt{-1}} 

\newcommand{\Spec}{\operatorname{Spec}} 
\newcommand{\Jac}{\operatorname{Jac}} 
\newcommand{\Lie}{\operatorname{Lie}}
\newcommand{\Hom}{\operatorname{Hom}} 
\newcommand{\Aut}{\operatorname{Aut}} 
\newcommand{\pr}{\operatorname{pr}} 
\newcommand{\GM}{\operatorname{GM}} 
\newcommand{\QDM}{\operatorname{QDM}} 
\newcommand{\ch}{\operatorname{ch}} 
\newcommand{\End}{\operatorname{End}}
\newcommand{\Res}{\operatorname{Res}}

\def\parfrac#1#2{\frac{\partial #1}{\partial #2}}
\def\corr#1{\left\langle #1 \right\rangle} 
\def\pair#1#2{\langle #1, #2 \rangle} 

\begin{document} 
\title{Quantum D-modules of toric varieties and oscillatory integrals}
\author{Hiroshi Iritani}
\email{iritani@math.kyoto-u.ac.jp} 
\address{Department of Mathematics, Graduate School of Science, 
Kyoto University, Kitashirakawa-Oiwake-cho, Sakyo-ku, 
Kyoto, 606-8502, Japan}
\maketitle

\begin{abstract} We review mirror symmetry for the quantum cohomology 
D-module of a compact weak-Fano toric manifold. 
We also discuss the relationship to the GKZ system, the Stanley-Reisner ring, 
the Mellin-Barnes integrals, and the $\hGamma$-integral structure.  
\end{abstract} 

\section{Introduction} 
The purpose of the present notes is to give a concise review on  
mirror symmetry for the quantum D-modules of toric varieties, 
as proposed by Givental \cite{Givental:ICM}. 
Our goal will be very modest: we restrict to weak-Fano compact 
toric manifolds and describe their mirror symmetry concretely 
and explicitly. 
We will also discuss the $\hGamma$-integral structure of  
the quantum D-modules and its role in mirror symmetry. 

A mirror of a Fano manifold $X$ is given by a Landau-Ginzburg 
model, that is, a complex manifold $Y$ equipped with a holomorphic function 
(potential) $W \colon Y \to \C$ (see \cite{Givental:ICM, 
Auroux:anticanonical, Golyshev:classification, 
KKP:Bogomolov, CCGGK}). 
Under mirror symmetry, it is expected that 
the quantum cohomology of $X$ is isomorphic to the Jacobi ring of $W$, 
and that 
the quantum cohomology D-module of $X$ 
is isomorphic to the mirror D-module whose solutions 
are oscillatory integrals associated with $W$ 
(see Table \ref{tab:mirror}). 
For a (weak-)Fano toric manifold $X$, the mirror $W$ is given by 
a family of Laurent polynomials whose Newton polytopes are 
given by the fan diagram of $X$, and the mirror correspondence 
as in Table \ref{tab:mirror} is by now well-established. 
We will explain the relationship between the mirror D-module 
and the Gelfand-Kapranov-Zelevinsky (GKZ) system 
and how the (quantum) Stanley-Reisner ring arises 
as a limit of the mirror D-module. 
The GKZ system or the Stanley-Reisner rings can be regarded 
as intermediate objects between quantum cohomology of 
toric varieties and their mirrors.

\renewcommand{\arraystretch}{1.4}
\begin{table}[h] 
\caption{Mirror correspondence.}
\label{tab:mirror}
\begin{tabular}{c|c| c} 
Fano manifold  & Landau-Ginzburg  & toric case \\ 
$X$ & model $W\colon Y\to \C$ & 
\\ \hline \hline 
quantum cohomology & Jacobi ring & quantum Stanley-\\ 
ring & (or $\mathbb{H}(Y,(\Omega_Y^\bullet,dW))$) & Reisner ring  
\\ \hline 
quantum D-module  & D-module of oscillatory  
& GKZ system 
\\ 
$z Q\parfrac{}{Q} + p\star$ & integrals 
$\int e^{W/z} \Omega$ & 
$\Box_d I = 0$ \\
\end{tabular}
\end{table} 

 Even for Fano manifolds other than toric varieties, 
the mirror space $Y$ is often given by 
(a partial compactification of) the algebraic torus 
$(\C^\times)^n$ and 
the potential $W$ is given by a Laurent polynomial on it; 
however, the coefficients of the mirror Laurent polynomial $W$ 
are very special (see e.g.~\cite{Rietsch:GmodP, CCGK}). 
Many people have observed (see \cite{BCFKvS}) that, 
when a Fano manifold $X$ admits a flat degeneration 
to a (possibly singular) toric variety $X'$, the mirror of $X$ is often 
given as a special subfamily of the Laurent polynomials that 
are mirror to a crepant resolution of $X'$. 
Such a subfamily corresponds to a certain stratum 
in the discriminant loci of the GKZ system for $X'$. 
Therefore, toric mirror symmetry can be viewed as 
a generic phenomenon whereas mirror symmetry 
for more general Fano manifolds can be viewed as its 
specialization. 
In this article, we do not delve into mirror symmetry for 
non-toric Fano manifolds; we refer the reader 
to the article of Golyshev \cite{Golyshev:techniques} 
in the same volume for the differential equations aspects 
of the story\footnote{This article was originally written 
as an appendix to his article, but became independent 
in the end.}. 

In the last two sections, we will discuss the $\hGamma$-integral 
structure \cite{Iritani:Integral, KKP:Hodge} 
in the context of toric mirror symmetry. 
The $\hGamma$-integral structure is a certain integral local system 
underlying the quantum D-module, which is defined by the $\hGamma$-class 
and the topological $K$-group. Conjecturally this corresponds to the 
natural integral structure on the mirror.  
The $\hGamma$-class involves transcendental numbers such as the 
Riemann $\zeta$-values and its origin is quite 
mysterious from a viewpoint of curve counting or symplectic topology. 
On the other hand, the existence of such a structure had been suggested 
since the beginning of mirror symmetry. Candelas-de la Ossa-Green-Parkes 
\cite{CdOGP} already observed that $\chi(X) \zeta(3)$  
appears in the asymptotics of periods of mirror Calabi-Yau 
threefolds near the large complex structure limit (LCSL). 
Hosono-Klemm-Theisen-Yau \cite{HKTY:CICY} observed 
``remarkable identities'' that relate certain characteristic numbers 
of complete intersection Calabi-Yau manifolds with 
hypergeometric solutions of the Picard-Fuchs equation of 
the mirror family. 
This observation led Libgober \cite{Libgober} 
to introduce the (inverse) $\hGamma$-class of a complex manifold. 
Later, Hosono \cite{Hosono:centralcharges} made this connection 
more precise in terms of central charges and homological mirror symmetry. 

When a Fano manifold $X$ is mirror to the Landau-Ginzburg model 
$W\colon Y \to \C$, 
the compatibility of the $\hGamma$-integral 
structure and mirror symmetry would imply the following 
``mirror symmetric Gamma conjecture'' (cf.~\cite{GGI,Galkin-Iritani}, 
see \S\ref{sec:Gamma})
\[
\int_{\Gamma \subset Y} 
e^{-t W} \Omega \sim \int_X t^{-c_1(X)} \cdot \hGamma_X 
\quad 
\text{as $t\to 0$.} 
\]
This says that the asymptotics of the mirror oscillatory integral  
is described in terms of the $\hGamma$-class of $X$. 
The integration cycle $\Gamma$ here should be a Lagrangian 
section of the Strominger-Yau-Zaslow fibration \cite{SYZ} 
(i.e.~mirror to the structure sheaf of $X$). 
In the last section, we will introduce Mellin-Barnes 
integral representations for the mirror oscillatory integrals 
(in the toric case) 
and explain these $\hGamma$-phenomena for toric varieties. 

\section*{Acknowledgements} 
I thank Vasily Golyshev for many helpful discussions. 
This work is supported by JSPS Kakenhi Grant Number 
16K05127, 16H06335, 16H06337, and 17H06127. 
The essential part of this work was done while I was in residence 
at the Hausdorff Research Institute for Mathematics at Bonn 
in January 2018, and at the Mathematical 
Sciences Research Institute at Berkeley, California in April 2018. 
The latter stay was supported by the National Science Foundation 
under Grant No.~DMS-1440140. I thank both institutes for 
providing excellent working conditions. 

\section{Preliminaries on toric varieties} 
A toric variety is a GIT quotient of a vector space 
$\C^m$ by a torus $K \cong (\C^\times)^k$, where the torus 
$K$ acts on $\C^m$ via an injective group homomorphism 
$K \to (\C^\times)^m$. 
Let $D_1,\dots,D_m \in \Hom(K,\C^\times)$ denote the characters 
defining the homomorphism $K \to (\C^\times)^m$. 
A GIT quotient is given by choosing a stability condition 
$\omega \in \frk_\R^*$, where $\frk_\R$ denotes the Lie algebra 
of the maximally compact subgroup $K_\R$ of $K$. 
We assume that 
\begin{itemize} 
\item[(a)] $\omega$ lies in 
the cone $\sum_{i=1}^m \R_{\ge 0} D_i$; 
\item[(b)] if $\omega$ lies in $\sum_{i\in I} \R_{>0} D_i$ 
for a subset $I\subset \{1,\dots,m\}$, then 
$\{D_i\}_{i\in I}$ 
spans $\frk_\R^*$ over $\R$; 
\item[(c)] the cone $\sum_{i=1}^m \R_{\ge 0} D_i$ is strictly convex.
\end{itemize}  
Then the toric variety corresponding to $\omega$ is defined by 
\[
X_\omega := U_\omega/K
\] 
where $U_\omega:=\C^m \setminus \bigcup_I \C^I$ 
with $I$ ranging over subsets of $\{1,\dots,m\}$ 
such that $\omega \notin \sum_{i\in I} \R_{>0} D_i$.  
The conditions (a), (b), (c) respectively ensure that 
$X_\omega$ is non-empty, 
that $X_\omega$ has at worst orbifold singularities, 
and that $X_\omega$ is compact. 
The space $\frk_\R^*$ of stability conditions has 
a fan structure, called the GKZ fan. A maximal cone of the GKZ fan 
is given by the closure of a connected component of 
the set of $\omega \in \frk_\R^*$ 
satisfying the conditions (a)--(c). The toric variety $X_\omega$ 
depends only on the maximal cone $A_\omega$ to which $\omega$ belongs. 

We assume that $X_\omega$ is weak-Fano, or equivalently, 
that $c_1(X_\omega)$ is nef. 
For simplicity of notation and exposition, 
we will also assume that $X_\omega$ is a smooth manifold 
(without orbifold singularities) and that 
$\{z_i=0\} \cap U_\omega$ is non-empty for all $i=1,\dots,m$, 
where $z_i$ is the $i$th co-ordinate on $\C^m$. 
With these assumptions, we have $\frk_\R^* \cong H^2(X_\omega,\R)$ 
and the closure of the ample cone of $X_\omega$ corresponds to the 
maximal cone $A_\omega$ of the GKZ fan. 
Under this isomorphism, the character $D_i$ corresponds to 
the class of the toric divisor 
$\{z_i = 0\}$ in $X_\omega$. 
\begin{remark}
Although we restrict to smooth toric manifolds in this article, 
all the results discussed in this paper can be extended to 
toric orbifolds
\cites{Iritani:Integral, Mann-Reichelt, CCLT, CCIT:MT, CCIT:MS}. 
In the orbifold case, 
it is important to allow $\{z_i=0\} \cap U_\omega = \emptyset$ 
for some $i$ as such indices correspond to twisted sectors.   
\end{remark} 

A toric variety can be also described in terms of a fan. 
Let $\frk_\Z= \Hom(\C^\times, K)$ denote the cocharacter lattice of $K$.  
Consider the natural map $\frk_\Z \to \Z^m$ induced 
by the inclusion $K\to (\C^\times)^m$ and complete it to 
a short exact sequence 
\[
0\to \frk_\Z \to \Z^m \to N \to 0
\]
where $N := \Z^m/\frk_\Z$ is a free abelian group of rank $n = m - k$. 
The fan $\Sigma_\omega$ of $X_\omega$ is defined on the vector space 
$N_\R= N\otimes \R$: the image $b_i \in N$ of $e_i \in \Z^m$ 
gives a generator of a 1-dimensional cone of $\Sigma_\omega$, 
and the cone $\sum_{i\in I} \R_{\ge 0} b_i$ belongs to $\Sigma_\omega$ 
if and only if $\omega \in \sum_{i \notin I} \R_{>0} D_i$.

\section{Quantum D-modules} 
The cocharacter lattice $\frk_\Z = \Hom(\C^\times, K)$ 
is identified with $H_2(X_\omega,\Z)$ 
and the dual cone $A_\omega^\vee\subset \frk_\R$ 
of $A_\omega$ is identified with the cone of curves. 
The (small) quantum product $\star$ is a commutative and 
associative product 
on the space $H^*(X_\omega) \otimes 
\C[\![A_\omega^\vee \cap \frk_\Z]\!]$ 
defined by 
\[
(\alpha\star \beta, \gamma) = \sum_{d\in A_\omega^\vee \cap \frk_\Z} 
\corr{\alpha,\beta,\gamma}_{0,3,d} Q^d \qquad 
\text{for all $\alpha,\beta,\gamma \in H^*(X_\omega)$},  
\] 
where $\corr{\alpha,\beta,\gamma}_{0,3,d}$ is the genus-zero, 
three-points, degree-$d$ Gromov-Witten invariant 
and $Q^d$ denotes 
the element of  $\C[\![A_\omega^\vee \cap \frk_\Z]\!]$ 
corresponding to $d$. 
The Dubrovin connection is a flat connection 
on the trivial $H^*(X_\omega)$-bundle over 
$\Spec \C[\![A_\omega^\vee \cap \frk_\Z]\!]$ given by 
\begin{equation} 
\label{eq:Dubrovin_conn} 
\nabla_{\xi Q\parfrac{}{Q}} = \xi Q\parfrac{}{Q} + 
\frac{1}{z}(\xi \star), \qquad 
\xi \in H^2(X_\omega) 
\end{equation} 
where $z$ is a formal parameter and $\xi Q\parfrac{}{Q}$ 
is a derivation of $\C[\![A_\omega^\vee \cap \frk_\Z]\!]$ 
defined by $(\xi Q\parfrac{}{Q}) Q^d = \pair{\xi}{d} Q^d$. 
This is a flat connection with logarithmic singularities. 
The Dubrovin connection $z \nabla$ (multiplied by $z$) acts on the space 
\[
\QDM(X_\omega) = H^*(X_\omega) \otimes 
\C[z][\![A^\vee_\omega \cap \frk_\Z]\!] 
\]
which we call the quantum D-module. 
It is not known in general whether the quantum product 
$\star$ converges or not. 
For toric varieties, it is known that 
the quantum product converges and hence 
the quantum D-module extends to an actual analytic 
neighbourhood of the origin ``$Q=0$'' in 
$\Spec \C[A_\omega^\vee \cap \frk_\Z]$. 
The point $Q=0$ is called the large radius limit point. 
\begin{remark} 
\label{rem:conn_z}
The Dubrovin connection can be also extended in the $z$-direction. 
The connection in the $z$-direction is given by 
\[
\nabla_{z\partial_z} = z\parfrac{}{z} - \frac{1}{z} c_1(X_\omega)\star 
+ \mu 
\]
where $\mu\in \End(H^*(X_\omega))$ is defined by 
$\mu(\alpha) = (p-\frac{n}{2}) \alpha$ 
for $\alpha \in H^{2p}(X_\omega)$ with $n= \dim_\C X_\omega$. 
\end{remark}

\section{Mirror D-modules}
\label{sec:mirror_D-mod} 
We have the exact sequence of tori 
$1 \to K \to (\C^\times)^m \to T \to 1$ 
where $T := (\C^\times)^m/K \cong N\otimes \C^\times$ 
is a torus acting on $X_\omega$ 
with an open dense orbit.  
Consider the exact sequence 
$1 \to \chT \to (\C^\times)^m \to \chK \to 1$ of dual tori. 
The mirror Landau-Ginzburg model of a toric variety $X_\omega$ 
is given by the family of tori 
$\pr\colon (\C^\times)^m \to \chK$ together with a 
potential function $W \colon (\C^\times)^m \to \C$ 
defined by $W = u_1+\cdots + u_m$, where $u_i$ denotes 
the $i$th co-ordinate on $(\C^\times)^m$. 
\[
\xymatrix{
(\C^\times)^m \ar[rrr]^{W=u_1+\cdots+u_m} \ar[d]^{\pr} 
&&& \C \\ 
\chK & && 
}
\]
Choosing a splitting of the sequence, we can also write 
$W_q=W|_{\pr^{-1}(q)}$ as 
\[
W_q=q^{l_1} x^{b_1} + \cdots + q^{l_m} x^{b_m}
\]
where $q \in \chK$, $x \in \chT$ and recall that $b_i$ are generators 
of 1-dimensional cones of the fan $\Sigma_\omega$. 
By varying $q\in \chK$, $W_q$ can represent any 
Laurent polynomial in $x\in \chT$ having $b_1,\dots,b_m$ as exponents. 
Hence the mirror of a toric variety can be thought of 
as generic Laurent polynomials. 

\begin{remark} 
We denote the B-model co-ordinates by $q$ and 
the A-model co-ordinates by $Q$.  
These co-ordinates are related by the mirror map $\psi$ below. 
\end{remark} 

Using the GKZ fan on $\frk_\R^*$ and 
its `preimage fan' on $\R^m$ (whose maximal cones 
are $(\R_{\ge 0})^m \cap \pi^{-1}(A)$, where 
$A$ is a maximal cone the GKZ fan and 
$\pi \colon \R^m \to \frk_\R^*$ is the natural map 
given by $D_1,\dots,D_m$), 
we can partially compactify 
the family $\pr \colon (\C^\times)^m \to \chK$ to a map between 
toric varieties $\pr\colon \cY \to \cM$ and 
$W$ extends to a regular function $\cY\to \C$. 
\[
\xymatrix{
(\C^\times)^m \ar[d] \ar@{^{(}->}[r] 
& \cY \ar[rr]^{W} \ar[d]^{\pr} && \C \\ 
\chK \ar@{^{(}->}[r] & \cM & & 
}
\]
The maximal cone $A_\omega$ defines a torus-fixed point 
$\bz_\omega \in \cM$ which is mirror to the large-radius limit 
point $Q=0$ of the quantum cohomology of $X_\omega$.  
Let $\cY_q$ denote the fibre of $q\in \cM$ and write $W_q = W|_{\cY_q}$. 
Givental introduced oscillatory integrals  
\begin{equation} 
\label{eq:oscint} 
\int_{\Gamma \subset \cY_q} e^{W_q/z} \Omega 
\end{equation} 
as mirrors of the quantum D-module, where $\Gamma$ is a 
non-compact Morse cycle for $\Re(W_q/z)$ 
and $\Omega$ is a holomorphic volume form 
on the fibre $\cY_q$. 
Introduce the log structures on $\cY$ and $\cM$ 
given by their toric boundaries and 
let $\Omega_{\cY/\cM}^\bullet$ denote the relative logarithmic 
de Rham complex. The integrands $\Omega$ of 
\eqref{eq:oscint} can be naturally viewed as elements of the twisted 
de Rham cohomology: 
\[
\GM(W) = 
\pr_*H^{\rm top}(\Omega_{\cY/\cM}^\bullet[z],zd + dW\wedge). 
\]
This is equipped with the (logarithmic) 
Gauss-Manin connection and the higher residue 
pairing; such structures were introduced 
by K.~Saito \cite{Saito:higher_residue} in singularity theory. 
The Gauss-Manin connection is a map 
\[
\nabla \colon \GM(W) \to \frac{1}{z} \GM(W) \otimes_{\cO_\cM} 
\Omega^1_\cM \oplus \GM(W) \frac{dz}{z^2}
\]
which has the same pole structure along $z=0$ as the Dubrovin connection 
(here $\Omega^1_\cM$ denotes the sheaf of logarithmic 1-forms). 
When we choose a splitting of the sequence $1\to \chT \to (\C^\times)^m 
\to \chK \to 1$ and choose co-ordinates $x=(x_1,\dots,x_n)$ on $\chT
\cong (\C^\times)^n$ and $q=(q_1,\dots,q_k)$ on 
$\chK\cong (\C^\times)^k$, the connection $\nabla$ is given by\footnote
{Here we twist the Gauss-Manin connection by $-\frac{n}{2}\frac{dz}{z}$ 
so that it is compatible with the Dubrovin connection.}
\[
\nabla [f \Omega_0] = \sum_{a=1}^k 
\left[\left(\partial_a f+ \frac{\partial_a W}{z}
 f \right)\Omega_0
\right] 
\frac{dq_a}{q_a} 
+ 
\left[\left(z\partial_z f - \frac{W}{z} f-\frac{n}{2} f \right) \Omega_0
\right] \frac{dz}{z}
\]
where $\partial_a = q_a\parfrac{}{q_a}$, $f=f(x,q,z) \in \cO_\cY[z]$ 
and $\Omega_0$ is the standard relative volume form 
of the family $\pr\colon \cY \to \cM$: 
\begin{equation} 
\label{eq:Omega_0} 
\Omega_0 = \frac{dx_1}{x_1} \wedge \cdots \wedge \frac{dx_n}{x_n}.  
\end{equation} 
The higher residue pairing is a map 
\[
(-)^*\GM(W) \otimes_{\cO_\cM[z]} \GM(W) \to \cO_\cM[\![z]\!] 
\]
which coincides with the residue pairing when restricted to $z=0$, 
where $(-)$ denotes the automorphism of $\cO_\cM[z]$ 
sending $f(q,z)$ to $f(q,-z)$. 
For a concrete description of the higher residue pairing, 
we refer the reader to \cite[\S 6]{CCIT:MS}. 
There exists a Zariski open subset $U$ of $\cM$ containing 
$\bz_\omega$ such that $\GM(W)$ is locally free and 
coherent (as an $\cO$-module) over $U\times \C_z$ 
\cite{Iritani:Integral, Reichelt-Sevenheck, Mann-Reichelt}. 

\begin{theorem}[\cites{Givental:ICM, Givental:mirror_thm_toric, 
Iritani:Integral, Douai-Mann, Reichelt-Sevenheck, Mann-Reichelt, CCIT:MS, 
Mochizuki:twistor_GKZ, Iritani:big_equivariant}] 
\label{thm:mirror_isom}
Over an analytic open neighbourhood of $\bz_\omega$, 
the mirror D-module $\GM(W)$ is isomorphic to 
the pull-back of the quantum D-module $\QDM(X_\omega)$ 
by a map $\psi \in \Aut(\C[\![A_\omega^\vee \cap \frk_\Z]\!])$, 
called the mirror map. 
Under this isomorphism, the higher residue pairing corresponds 
to the Poincar\'e pairing. 
\end{theorem} 
The mirror map $\psi$ above extends 
to a local isomorphism between analytic neighbourhoods of 
$\bz_\omega$, and therefore the pull-back by $\psi$ 
makes sense over an analytic neighbourhood. 
If $\cM$ contains two different 
large radius limit points $\bz_{\omega_1}$, $\bz_{\omega_2}$ 
such that $X_{\omega_1}, X_{\omega_2}$ are weak Fano, the 
theorem implies that the quantum D-modules of $X_{\omega_1}$ 
and $X_{\omega_2}$ are isomorphic under analytic continuation; this 
is an example where the crepant transformation conjecture holds 
\cite{CIT,CIJ}. 

\begin{remark} 
\label{rem:historical} 
The above theorem is essentially due to Givental \cite{Givental:ICM, 
Givental:mirror_thm_toric}. 
Givental introduced the mirror potential function $W$ 
(see also Hori-Vafa \cite{Hori-Vafa}) and 
expressed solutions of the quantum D-module in terms of oscillatory integrals. 
Givental's mirror theorem \cite{Givental:mirror_thm_toric} 
(see also Lian-Liu-Yau \cite{LLY:I, LLY:II})  
is stated as the equality between the $J$-function and the $I$-function, 
where the $J$-function is a solution to the quantum D-module 
and the $I$-function is a solution to the mirror D-module 
(see the next section). 
Mirror symmetry as an isomorphism of D-modules has been studied 
in details by 
\cite[Proposition 4.8]{Iritani:Integral}, 
\cite[Theorem 4.11]{Reichelt-Sevenheck}, 
\cite[Theorem 5.1.1]{Douai-Mann}, 
\cite[Theorem 7.43]{Mochizuki:twistor_GKZ}, 
\cite[Theorem 4.2]{Iritani:big_equivariant}, 
\cite[Theorem 6.4]{Mann-Reichelt}, 
\cite[Theorem 1.1]{CCIT:MS}. 
The logarithmic extension across the boundary divisors of $\cM$ 
has been studied by \cite{Reichelt-Sevenheck, Douai-Mann,
Mochizuki:twistor_GKZ, Mann-Reichelt} in terms of (GKZ) D-modules; 
the logarithmic extension of the mirror Landau-Ginzburg model itself 
was discussed in \cite{Mann-deGregorio, Iritani:big_equivariant, CCIT:MS}. 
In this paper, we restrict to the small quantum cohomology, but 
we also have mirror symmetry for the big quantum cohomology,  
see \cite{Barannikov:projective, Douai-Sabbah:I, Douai-Sabbah:II, 
Iritani:big_equivariant, CCIT:MS}. 
\end{remark} 

\begin{example} 
\label{exa:Pn}
The mirror family of $\PP^n$ is given by the diagram 
\[
\xymatrix{
\C^{n+1} \ar[d]^{\pr} 
\ar[rrr]^{W=u_0+\dots+u_n} && & \C \\ 
\C && & 
} 
\] 
with $q= \pr(u_0,\dots,u_n) = u_0u_1\cdots u_n$, 
where $\cY = \C^{n+1}$ and $\cM = \C$. 
We can write $W_q=W|_{\pr^{-1}(q)}$ as $W_q = x_1+\cdots+x_n+ 
\frac{q}{x_1\cdots x_n}$ by setting $x_i := u_i$ 
for $1\le i\le n$. 
\end{example} 

\section{The GKZ system and hypergeometric solutions} 
\label{sec:GKZ} 
The mirror D-module $\GM(W)$ can be described in terms of 
the Gelfand-Kapranov-Zelevinsky (GKZ) system \cite{GKZ}. 
To describe the GKZ system explicitly, we introduce a basis 
 $p_1,\dots,p_k$ of $\frk_\Z^*\cong H^2(X_\omega,\Z)$ 
and write $D_i = \sum_{a=1}^k m_{ia} p_a$. 
The elements $p_1,\dots,p_k$ define co-ordinates 
$q_1,\dots,q_k$ on $\chK \cong (\C^\times)^k$. 
We write $q^d = \prod_{a=1}^k q_a^{p_a \cdot d}$ for $d\in \frk_\Z$ 
and set $\partial_a = q_a \parfrac{}{q_a}$.  
Over the open torus $\chK \subset \cM$, the mirror D-module 
$\GM(W)$ is generated by the standard relative  
volume form $\Omega_0$ (see \eqref{eq:Omega_0}) 
as a D-module, or more precisely, as a module over the ring 
$\cO_{\chK}[z]\langle z\partial_1,\dots,z\partial_k\rangle$, 
where $z\partial_a$ acts by the Gauss-Manin connection 
$z \nabla_{\partial_a}$. 
All the D-module relations of $\Omega_0$ 
are generated by $\Box_d [\Omega_0]=0$ with 
$d\in \frk_\Z$, where 
\begin{equation} 
\label{eq:GKZ_op} 
\Box_d := \prod_{i: D_i\cdot d>0} \prod_{j=0}^{D_i \cdot d-1} 
\left(\sum_{a=1}^k m_{ia} z\partial_a - jz\right) - q^d 
\prod_{i: D_i \cdot d<0} \prod_{j=0}^{-D_i \cdot d -1} 
\left(\sum_{a=1}^k m_{ia} z \partial_a - jz\right). 
\end{equation} 
This system of differential equations is called 
the GKZ system. 
Givental's $I$-function \cite{Givental:mirror_thm_toric} 
is a cohomology-valued function annihilated by $\Box_d$: 
\[
I(q,z) = \sum_{d\in A_\omega^\vee \cap \frk_\Z} 
q^{d+p/z} \prod_{i=1}^m \frac{\prod_{j=-\infty}^{0}(D_i+jz)}{
\prod_{j=-\infty}^{D_i\cdot d} (D_i + jz) } 
\]
where $q^{p/z} = \prod_{a=1}^k e^{p_a\log q_a/z}$. 
The components of the $I$-function 
form a basis of solutions to 
the GKZ system near $\bz_\omega$ 
\cite{Iritani:Integral, Reichelt-Sevenheck, Mann-Reichelt}. 
The mirror map $\psi$ is determined 
by the $z^{-1}$-expansion of the $I$-function 
\[
I(q,z) = 1 + \frac{1}{z} \sum_{a=1}^k p_a \log \psi_a(q) 
+ O(z^{-2}) 
\]
where we write the mirror map in the form 
$Q_a = \psi_a(q)$, $a=1,\dots,k$ 
using the A-model co-ordinates $Q_1,\dots,Q_k$ 
dual to $p_1,\dots,p_k$. 
We have $\psi_a(q) = q_a + \text{higher order terms}$; 
in the Fano case the mirror map is trivial $\psi_a(q) = q_a$. 
\begin{remark}[\cite{GKZ,Adolphson:hg,Iritani:Integral}] 
The rank of the GKZ system (at a generic point in $\chK$) 
is equal to the normalized volume of the fan polytope, 
that is, the convex hull of $b_1,\dots,b_m$ in $N_\R$; 
this is also equal to $\dim H^*(X_\omega)$. 
\end{remark} 

\begin{example}[continuation of Example \ref{exa:Pn}]
\label{exa:Pn_I}
In the case of $\PP^n$, the mirror oscillatory integral satisfies 
the differential equation: 
\[
\left(\left( z q\partial_q \right)^{n+1} - q\right) 
\int_\Gamma e^{\left(x_1+\cdots+x_n+\frac{q}{x_1\cdots x_n}\right)/z} 
\frac{dx_1\cdots dx_n}{x_1\cdots x_n} =0.  
\]
The $I$-function is given by 
\[
I(q,z) = \sum_{d=0}^\infty 
\frac{q^{p/z+d}}{\prod_{j=1}^d (p+jz)^{n+1}}
\]
(with $p= c_1(\cO(1))$) and the mirror map is trivial: $Q= q$. 
The relationship between these two solutions will be discussed in 
\S\ref{sec:Gamma}-\ref{sec:MB}, see Remark \ref{rem:I_oscillatory} 
and Example \ref{exa:Pn_Mellin}.  
\end{example} 

\section{The Jacobi ring and the Stanley-Reisner ring} 
A ring structure arises from these D-modules in the $z\to 0$ limit. 
The $z\to 0$ limit of the Dubrovin connection 
$z \nabla_{\xi Q\parfrac{}{Q}}$ 
is the quantum product by $\xi\in H^2(X_\omega)$ (see \eqref{eq:Dubrovin_conn}), 
and the quantum multiplication by divisors generate 
the quantum cohomology ring. 
On the other hand, the $z\to 0$ limit of the mirror D-module is 
given by the (logarithmic) Jacobi ring of $W$: 
\[
\Jac(W) = \pr_* \left(
\cO_{\cY}\Big/\left \langle \textstyle x_1\parfrac{W}{x_1}, 
\dots, x_n \parfrac{W}{x_n} \right\rangle\right)  
\] 
where we choose $\C^\times$-co-ordinates $(x_1,\dots,x_n)$ 
of $\chT \cong (\C^\times)^n$ and 
$x_i \parfrac{}{x_i}$ defines a relative tangent vector field 
of the family $\cY \to \cM$. 
Theorem \ref{thm:mirror_isom} induces an isomorphism 
\[
\psi^*QH^*(X_\omega) \cong \Jac(W) 
\]
between the quantum cohomology ring and the Jacobi ring. 
This isomorphism has been studied by many people, 
see \cite{Batyrev, McDuff-Tolman, 
FOOO:I, FOOO:II, FOOO:Ast, Gonzalez-Woodward:tmmp, 
Ritter,Smith} together with the references 
in Remark \ref{rem:historical}. 

Over the affine open chart $\cM_\omega :=
\Spec(\C[A_\omega^\vee \cap \frk_\Z])$ of $\cM$, 
the Jacobi ring is a quotient of the quantum Stanley-Reisner ring 
associated with the fan $\Sigma_\omega$. 
For $v\in N$, choose a cone $\sum_{i\in I} \R_{\ge 0} b_i$ 
of the fan $\Sigma_\omega$ containing $v$ and write 
$v = \sum_{i\in I} v_i b_i$. We set $v_i =0$ for $i\notin I$. 
Set $w_v := \prod_{i=1}^m u_i^{v_i}$. Then 
$H^0(\cM_\omega,\pr_*\cO_{\cY})$ is a free 
$\C[A_\omega^\vee \cap \frk_\Z]$-module 
generated by $w_v$ with $v\in N$. 
The product structure is given by 
\begin{equation} 
\label{eq:qSR}
w_{v} w_{v'} = q^{\ell(v,v')} w_{v+v'} 
\end{equation} 
where $\ell(v,v')\in A_\omega^\vee \cap \frk_\Z$ is the element 
given by the linear relation $(v_i + v'_i - (v+v')_i)_{i=1}^m$ 
among $b_i$'s via the exact sequence 
$0\to \frk_\Z \to \Z^m \to N \to 0$. 
This ring $H^0(\cM_\omega,\pr_*\cO_{\cY})$ is called the 
quantum Stanley-Reisner ring.  
The Jacobian ideal gives additional linear relations among 
$u_1,\dots,u_m$ given by: 
\begin{equation} 
\label{eq:linear_relation} 
\sum_{i=1}^m u_i b_i =0. 
\end{equation} 
The relations \eqref{eq:qSR}, \eqref{eq:linear_relation} 
define the Jacobi ring $\Jac(W)$; they are also known as 
the relations of Batyrev's quantum ring \cite{Batyrev}. 
At the large-radius limit $q\to \bz_\omega$, the quantum 
Stanley-Reisner relations \eqref{eq:qSR}  
reduce to 
\[
w_v w_{v'} = \begin{cases} 
w_{v+v'}  & 
\text{if $v$ and $v'$ lie in a common cone of $\Sigma_\omega$;} \\
0 & \text{otherwise.} 
\end{cases} 
\]
These relations define the Stanley-Reisner ring. 
As is well-known in toric geometry, the Stanley-Reisner ring modulo 
the linear relations \eqref{eq:linear_relation} 
is isomorphic to $H^*(X_\omega)$. 

\begin{remark} 
Over an analytic neighbourhood of $\bz_\omega$, 
$\pr_*\cO_{\cY}$ is isomorphic to the $T$-equivariant   
quantum cohomology of $X_\omega$ and the left-hand side of 
the linear relations \eqref{eq:linear_relation}
correspond to the $T$-equivariant parameters. 
Moreover, the multiplication by the co-ordinates $x_1,\dots,x_n$ of 
$\chT$ corresponds to 
the Seidel representation \cite{Seidel} on the equivariant 
quantum cohomology. 
McDuff-Tolman \cite{McDuff-Tolman} used the Seidel representation 
to determine a presentation of the quantum cohomology ring 
of a toric manifold. 
In \cite{Iritani:big_equivariant}, we showed how equivariant 
mirror symmetry for toric manifolds 
follows almost tautologically from the Seidel representation 
and shift operators \cite{Braverman-Maulik-Okounkov}. 
\end{remark} 

\section{$\hGamma$-integral structure}
\label{sec:Gamma} 
The $\hGamma$-class \cite{Libgober} 
of an almost complex manifold $X$ is defined to 
be the characteristic class 
\[
\hGamma_X = \prod_{i=1}^n \Gamma(1+\delta_i) 
\]
where $\delta_1,\dots,\delta_n$ are the Chern roots of 
the tangent bundle (so that $c(TX) = \prod_{i=1}^n (1+\delta_i)$) 
and $\Gamma(x)$ is Euler's $\Gamma$-function. 
The $\Gamma$-function $\Gamma(1+\delta)$ here  
should be expanded in Taylor series at $\delta=0$. 
For a toric variety $X = X_\omega$, it is given by  
\[
\hGamma_{X_\omega} = \prod_{j=1}^m 
\Gamma(1+D_j)
=\prod_{j=1}^m \exp\left(-\gamma D_j +\sum_{k=2}^\infty 
(-1)^k \frac{\zeta(k)}{k} D_j^k \right)  
\]
where $\gamma$ is the Euler constant, $\zeta(s)$ is the Riemann 
$\zeta$-function and 
recall that $D_j$ is the class of a prime toric divisor. 

The $\hGamma$-integral structure \cite{Iritani:Integral, 
KKP:Hodge} is an integral lattice in the 
space of (multi-valued) flat sections of the Dubrovin connection.  
For an element $E\in K^0(X_\omega)$ of the topological $K$-group, 
there exists a unique flat section $s_E(Q,z)$ of the Dubrovin connection 
(which is flat also in the $z$-direction, see Remark \ref{rem:conn_z}) 
such that 
\[
s_E(Q,z) \sim 
\frac{1}{(2\pi)^{n/2}}Q^{-p/z} z^{-\mu} z^{c_1(X_\omega)} 
\left( 
\hGamma_{X_\omega} \cup (2\pi\iu)^{\deg/2} \ch(E) \right) 
\]
as $Q \to 0$, where $Q^{-p/z}=\prod_{a=1}^k e^{-p_a \log Q_a/z}$ 
(we choose co-ordinates $Q_1,\dots,Q_k$ dual to 
$p_1,\dots,p_k$ as in \S\ref{sec:GKZ}) and $\mu$ is as in Remark 
\ref{rem:conn_z}. 
These flat sections $s_E(Q,z)$ span the $\hGamma$-integral structure. 
Its important properties are as follows: 
\begin{itemize} 
\item[(i)] this lattice is invariant under local monodromy 
around the large radius limit point; therefore it defines a $\Z$-local system 
underlying the quantum D-module; 
\item[(ii)] the Poincar\'e pairing of these flat sections coincide with 
the Euler pairing on the derived category: 
$(s_E(Q,e^{-\pi\iu}z),s_F(Q,z)) = \chi(E,F)$. 
\end{itemize} 
The second property follows from the identity 
$\Gamma(1+x) \Gamma(1-x) = \pi x/\sin(\pi x)$ 
(which implies that the $\hGamma_X$ is the ``half'' 
of the Todd class) and the Hirzebruch-Riemann-Roch formula. 

The mirror D-module $\GM(W)$ from \S\ref{sec:mirror_D-mod} 
has a natural integral structure\footnote{To be more precise, 
we twist the local system of relative cohomology 
by $(-2\pi z)^{-n/2}$ so that it is compatible 
with the Dubrovin connection.} 
dual to the relative homology 
$H_n(\cY_q,\{\Re(W_q/z)\ll 0\};\Z)$ 
via the oscillatory integral \eqref{eq:oscint}. 
Here note that an element of the relative homology gives an integration 
cycle in \eqref{eq:oscint}. 
These two integral structures coincide under mirror symmetry. 
\begin{theorem}[\cite{Iritani:Integral}]  
The $\hGamma$-integral structure coincides with the natural integral 
structure on the mirror D-module under the mirror isomorphism in 
Theorem $\ref{thm:mirror_isom}$.   
\end{theorem} 

This theorem follows from the following identity of periods  
of both sides: 
\begin{equation} 
\label{eq:equality_periods} 
 \int_{X_\omega} s_E(\psi(q),z) 
= \frac{1}{(2\pi z)^{n/2}}\int_{\Gamma(E)} e^{-W_q/z} \Omega_0 
\end{equation} 
where $E\mapsto \Gamma(E)$ is an isomorphism  
between $K^0(X_\omega)$ and $H_n(\cY_q, \{\Re(W_q/z)\gg 0\};\Z)$; 
this correspondence should be a shadow of homological 
mirror symmetry, i.e.~an equivalence between 
the derived category of coherent sheaves on $X_\omega$ 
and the Fukaya-Seidel category of $W_q$ (see \cite{Fang}): 
\[
D^b_{\rm coh}(X_\omega) \cong FS(Y,W_q).
\]
When $E$ is the structure sheaf $\cO$ on $X_\omega$, 
the corresponding cycle $\Gamma(\cO)$ is given by the positive real 
locus $(\R_{>0})^n$ in $\cY_q\cong \chT\cong (\C^\times)^n$ 
(when $q$ lies in the positive real locus $\Hom(\frk_\Z,\R_{>0})$ 
of $\chK = \Hom(\frk_\Z,\C^\times)$ and $z>0$). 
In this case, \eqref{eq:equality_periods} yields the following asymptotics: 
\begin{align}
\label{eq:asymp_Gamma}  
\int_{(\R_{>0})^n}e^{-W_q/z} \Omega_0 
& \sim \int_{X_\omega} q^{-p} z^{c_1(X_\omega)} \cup 
\hGamma_{X_\omega} 
\end{align} 
as $q\to \bz_\omega$ in the positive real locus and for $z>0$. 
Namely, the $\hGamma$-class appears in the $q\to \bz_\omega$ 
asymptotics of the exponential period of the mirror. 
It would be very interesting to study if such asymptotics hold for more 
general Fano manifolds and their mirrors. 

\begin{remark} 
\label{rem:I_oscillatory} 
The equality \eqref{eq:equality_periods} 
can be restated in terms of the $I$-function as follows: 
\[
\int_{X_\omega} 
\left(z^{c_1(X_\omega)} z^{\frac{\deg}{2}} I(q,-z)\right) 
\cup \hGamma_{X_\omega} (2\pi\iu)^{\frac{\deg}{2}}\ch(E) 
= \int_{\Gamma(E)} e^{-W_q/z} \Omega_0. 
\]
This formula expresses the oscillatory integral as an explicit  
linear combination of components of the $I$-function. 
\end{remark} 


\begin{remark} 
There is another version of ``Gamma conjecture'' for quantum 
cohomology of Fano manifolds formulated by 
Galkin, Golyshev and the author \cite{GGI}, 
which does not involve mirror symmetry. 
The quantum D-module of a Fano manifold has irregular singularities 
at ``$Q=\infty$'', and the conjecture is about the Stokes structure 
at $Q=\infty$ and the connection of solutions between $Q=0$ and $Q=\infty$. 
This version of the conjecture, if true, implies that the $\hGamma$-class 
can be recovered from the quantum cohomology of a Fano manifold. 
See also \cite{Galkin-Iritani, Golyshev-Zagier}. 
\end{remark} 

\begin{remark} 
Abouzaid, Ganatra, Sheridan and the author \cite{AGIS} 
recently proposed an approach to proving the 
asymptotics \eqref{eq:asymp_Gamma} 
(for Calabi-Yau mirror pairs) 
using Strominger-Yau-Zaslow 
picture \cite{SYZ} and tropical geometry. 
\end{remark} 

\section{Mellin-Barnes integral representation}
\label{sec:MB}
The oscillatory integrals \eqref{eq:oscint} give an integral 
representation of solutions to the GKZ system. 
There is another integral representation, Mellin-Barnes integral 
representation, which is ``Gale dual'' to \eqref{eq:oscint}. 
We shall regard $p_1,\dots,p_k \in \frk_\Z^*$ 
as co-ordinates on $\frk_\C = \Lie(K)$ which 
are Mellin-dual to $q_1,\dots,q_k$.  
Under the Mellin transformation $I(q) \mapsto 
\hI(p) = \int  I(q)q^{p} \frac{dq}{q}$, the differential operator 
$\partial_a =q_a\parfrac{}{q_a}$ 
corresponds to the multiplication by $-p_a$ and 
the multiplication by 
$q_a$ corresponds to the shift operator $T_a \colon 
p_b \mapsto p_b+\delta_{a,b}$. 
Thus the GKZ equations  $\Box_d I(q,z)=0$ (see \eqref{eq:GKZ_op}) 
with $d\in \frk_\Z$ are transformed into the difference equations: 
\[
\left( \prod_{i: D_i\cdot d>0} 
\prod_{j=0}^{D_i\cdot d -1} (-D_i - j)z  
- T^d \prod_{i: D_i \cdot d<0} 
\prod_{j=0}^{D_i \cdot d -1}
(-D_i-j) z \right) \hI(p,z) =0 
\]
where $D_i = \sum_{a=1}^k m_{ia} p_a\in \frk_\Z^*$ is
regarded as the pull-back of the standard co-ordinates on $\C^m$ 
via the inclusion $\frk_\C \hookrightarrow \C^m$ 
and $T^d = \prod_{a=1}^k T_a^{p_a\cdot d}$. 
It is easy to check that 
this system has the following simple solution:  
\[
\hI(p,z) = \prod_{i=1}^m (-z)^{D_i}\Gamma(D_i). 
\]
In fact, $\hI$ is (the restriction of) the Mellin transform of  
$e^{W/z} = \prod_{i=1}^m e^{u_i/z}$; 
here we recall 
$\int_0^\infty e^{u/z} u^{D} \frac{du}{u} =
(-z)^D \Gamma(D)$ with $z<0$ and $D>0$. 
By the inverse Mellin transformation, 
we get a solution to the GKZ system: 
\begin{equation}
\label{eq:MB}  
\frac{1}{(2\pi\iu)^k} 
\int_{C\subset \frk_\C} q^{-p} 
\left(\prod_{i=1}^m (-z)^{D_i} \Gamma(D_i) \right) 
dp_1\cdots dp_k 
\end{equation} 
where $C \subset \frk_\C$ is a suitable (non-compact) $k$-cycle 
so that the integral converges. 
For a suitable choice of $C$, \eqref{eq:MB} should coincide 
with an oscillatory integral $\int_{\Gamma} e^{W_q/z} \Omega_0$ 
(see Figure \ref{fig:Mellin} below), 
but the author does not know a precise choice of cycles 
in general.  
Via the residue calculation, such a formula would 
explain the $\hGamma$-class appearing in the leading 
asymptotics \eqref{eq:asymp_Gamma}, see 
Example \ref{exa:Pn_Mellin} below. 
\begin{figure}[h] 
\[
\xymatrix{
e^{W/z} = \prod_{i=1}^m e^{u_i/z} \ \text{on $(\C^\times)^m$} 
\ar[rr]^{\phantom{ABCD}\pr_*} 
\ar@{<->}[d]^{\text{Mellin}} 
&&  \int e^{W/z} \Omega \ \text{on $\chK$}  
\ar@{<->}[d]^{\text{Mellin}} 
\\ 
\prod_{i=1}^m (-z)^{D_i} \Gamma(D_i) \ \text{on $\C^m$} 
\ar[rr]^{\phantom{ABCDE}\text{restriction}}& & 
 \text{$\hI(p,z)$ on $\frk_\C$}  
}
\]
\caption{Oscillatory integral and its Mellin transform}
\label{fig:Mellin} 
\end{figure} 
\begin{example}[continuation of Example \ref{exa:Pn_I}] 
\label{exa:Pn_Mellin}
We consider the mirror oscillatory integral of $\PP^n$ again. 
The following method is borrowed from \cite{KKP:Hodge}. 
The Mellin transform of the oscillatory integral 
$\cI(q) = \int_{(\R_{>0})^n} 
e^{-(x_1+\cdots+x_n+\frac{q}{x_1\cdots x_n})/z} \frac{dq}{q}$ 
(with $q,z>0$) gives 
\begin{align*} 
\hcI(p) & = \int_0^\infty q^p \cI(q) \frac{dq}{q} \\ 
& = 
\int_{(\R_{>0})^{n+1}} (u_0\cdots u_n)^p 
e^{-(u_0+\cdots+u_n)/z} \frac{du_0}{u_0}\wedge 
\cdots \wedge \frac{du_n}{u_n} \\ 
& = z^{(n+1)p} \Gamma(p)^{n+1}. 
\end{align*} 
This coincides with $\hI(p,-z)$ as expected. 
Then the Mellin inversion formula gives 
\[
\cI(q) = \frac{1}{2\pi\iu} \int_{c-\iu \infty}^{c+\iu \infty} 
q^{-p} \hcI(p) dp 
\]
with $c>0$. By closing the integration contour to the left, 
we can express the right-hand side as the sum over 
residues at $p=0,-1,-2,\dots$, arriving at 
the asymptotics in \eqref{eq:asymp_Gamma}. 
\begin{align*} 
\cI(q) & = \sum_{d=0}^\infty \Res_{p=-d} 
\left( q^{-p}z^{(n+1)p}
\Gamma(p)^{n+1}  dp \right) \\ 
& \sim \Res_{p=0} 
\left(q^{-p} z^{(n+1)p} \Gamma(1+p)^{n+1} 
\frac{dp}{p^{n+1}}\right) 
= \int_{\PP^n} q^{-p} z^{c_1(\PP^n)} \cup 
\hGamma_{\PP^n}. 
\end{align*} 
\end{example} 

\begin{remark} 
The Mellin-Barnes integral representations \eqref{eq:MB} 
appear in 
physics literature as hemisphere partition functions 
(studied for more general gauged linear sigma models), 
see, e.g.~\cite{Hori-Romo}.  
\end{remark}

\bibliographystyle{amsplain} 

\begin{bibdiv}
\begin{biblist}  

\bib{AGIS}{article}{
author={Abouzaid, Mohammed}, 
author={Ganatra, Sheel}, 
author={Iritani, Hiroshi}, 
author={Sheridan, Nick}, 
title={The Gamma and Strominger-Yau-Zaslow conjectures: a tropical approach to periods}, 
note={\arxiv{1809.02177}}
}

\bib{Adolphson:hg}{article}{
   author={Adolphson, Alan},
   title={Hypergeometric functions and rings generated by monomials},
   journal={Duke Math. J.},
   volume={73},
   date={1994},
   number={2},
   pages={269--290},
   issn={0012-7094},
   doi={10.1215/S0012-7094-94-07313-4},
}

\bib{Auroux:anticanonical}{article}{
   author={Auroux, Denis},
   title={Mirror symmetry and $T$-duality in the complement of an
   anticanonical divisor},
   journal={J. G\"{o}kova Geom. Topol. GGT},
   volume={1},
   date={2007},
   pages={51--91},
   issn={1935-2565},
}

\bib{Barannikov:projective}{article}{
author={Barannikov, Serguei}, 
title={Semi-infinite Hodge structure and mirror symmetry 
for projective spaces}, 
note={\arxiv{math.AG/0010157}}, 
year={2001}
}

\bib{Batyrev}{article}{
   author={Batyrev, Victor V.},
   title={Quantum cohomology rings of toric manifolds},
   note={Journ\'{e}es de G\'{e}om\'{e}trie Alg\'{e}brique d'Orsay (Orsay, 1992)},
   journal={Ast\'{e}risque},
   number={218},
   date={1993},
   pages={9--34},
   issn={0303-1179},
}

\bib{BCFKvS}{article}{
   author={Batyrev, Victor V.},
   author={Ciocan-Fontanine, Ionu\c{t}},
   author={Kim, Bumsig},
   author={van Straten, Duco},
   title={Conifold transitions and mirror symmetry for Calabi-Yau complete
   intersections in Grassmannians},
   journal={Nuclear Phys. B},
   volume={514},
   date={1998},
   number={3},
   pages={640--666},
   issn={0550-3213},
   doi={10.1016/S0550-3213(98)00020-0},
}


\bib{Braverman-Maulik-Okounkov}{article}{
author={Braverman, Alexander},  
author={Maulik, Davesh}, 
author={Okounkov, Andrei}, 
title={Quantum cohomology of the Springer resolution}, 
journal={Adv. Math.}, 
volume={227}, 
number={1}, 
pages={421--458}, 
year={2011}
}

\bib{CdOGP}{article}{
   author={Candelas, Philip},
   author={de la Ossa, Xenia C.},
   author={Green, Paul S.},
   author={Parkes, Linda},
   title={A pair of Calabi-Yau manifolds as an exactly soluble
   superconformal theory},
   journal={Nuclear Phys. B},
   volume={359},
   date={1991},
   number={1},
   pages={21--74},
   issn={0550-3213},
   doi={10.1016/0550-3213(91)90292-6},
}

\bib{CIJ}{article}{
author={Coates, Tom}, 
author={Iritani, Hiroshi}, 
author={Jiang, Yunfeng}, 
title={The crepant transformation conjecture for toric complete
   intersections},
   journal={Adv. Math.},
   volume={329},
   date={2018},
   pages={1002--1087},
   issn={0001-8708},
   doi={10.1016/j.aim.2017.11.017},
}


\bib{CCIT:MT}{article}{
author={Coates, Tom}, 
author={Corti, Alessio}, 
author={Iritani, Hiroshi}, 
author={Tseng, Hsian-Hua}, 
title ={A mirror theorem for toric stacks}, 
journal={Compos. Math.}, 
volume={151}, 
number={10}, 
pages={1878--1912}, 
year={2015} 
}
\bib{CCIT:MS}{article}{
author ={Coates, Tom}, 
author ={Corti, Alessio}, 
author ={Iritani, Hiroshi}, 
author ={Tseng, Hsian-Hua}, 
title ={Hodge-theoretic mirror symmetry for toric stacks}, 
note ={\arxiv{1606.07254}, to appear in Journal of Differential Geometry}, 
}

\bib{CCGGK}{article}{
   author={Coates, Tom},
   author={Corti, Alessio},
   author={Galkin, Sergey},
   author={Golyshev, Vasily},
   author={Kasprzyk, Alexander},
   title={Mirror symmetry and Fano manifolds},
   conference={
      title={European Congress of Mathematics},
   },
   book={
      publisher={Eur. Math. Soc., Z\"{u}rich},
   },
   date={2013},
   pages={285--300},
}

\bib{CCGK}{article}{
   author={Coates, Tom},
   author={Corti, Alessio},
   author={Galkin, Sergey},
   author={Kasprzyk, Alexander},
   title={Quantum periods for 3-dimensional Fano manifolds},
   journal={Geom. Topol.},
   volume={20},
   date={2016},
   number={1},
   pages={103--256},
   issn={1465-3060},
   doi={10.2140/gt.2016.20.103},
}

\bib{CCLT}{article}{ 
author={Coates, Tom}, 
author={Corti, Alessio}, 
author={Lee, Yuan-Pin}, 
author={Tseng, Hsian-Hua}, 
title={The quantum orbifold cohomology of weighted projective spaces}, 
journal={Acta Math.}, 
volume={202},  
number={2}, 
pages={139--193}, 
year={2009}
} 

\bib{CIT}{article}{
   author={Coates, Tom},
   author={Iritani, Hiroshi},
   author={Tseng, Hsian-Hua},
   title={Wall-crossings in toric Gromov-Witten theory. I. Crepant examples},
   journal={Geom. Topol.},
   volume={13},
   date={2009},
   number={5},
   pages={2675--2744},
   issn={1465-3060},
   review={\MR{2529944}},
   doi={10.2140/gt.2009.13.2675},
}

\bib{Douai-Mann}{article}{
author={Douai, Antoine},  
author={Mann, Etienne}, 
title={The small quantum cohomology of a weighted projective space, 
a mirror D-module and their classical limits}, 
journal={Geom. Dedicata}, 
volume={164}, 
pages={187--226}, 
year={2013}
}

\bib{Douai-Sabbah:I}{article}{
   author={Douai, A.},
   author={Sabbah, C.},
   title={Gauss-Manin systems, Brieskorn lattices and Frobenius structures.
   I},
   language={English, with English and French summaries},
   booktitle={Proceedings of the International Conference in Honor of
   Fr\'{e}d\'{e}ric Pham (Nice, 2002)},
   journal={Ann. Inst. Fourier (Grenoble)},
   volume={53},
   date={2003},
   number={4},
   pages={1055--1116},
   issn={0373-0956},
}

\bib{Douai-Sabbah:II}{article}{
   author={Douai, Antoine},
   author={Sabbah, Claude},
   title={Gauss-Manin systems, Brieskorn lattices and Frobenius structures.
   II},
   conference={
      title={Frobenius manifolds},
   },
   book={
      series={Aspects Math., E36},
      publisher={Friedr. Vieweg, Wiesbaden},
   },
   date={2004},
   pages={1--18},
}

\bib{Fang}{article}{
author={Fang, Bohan}, 
title={Central charges of T-dual branes for toric varieties}, 
note={\arxiv{1611.05153}}
}

\bib{FOOO:I}{article}{
   author={Fukaya, Kenji},
   author={Oh, Yong-Geun},
   author={Ohta, Hiroshi},
   author={Ono, Kaoru},
   title={Lagrangian Floer theory on compact toric manifolds. I},
   journal={Duke Math. J.},
   volume={151},
   date={2010},
   number={1},
   pages={23--174},
   issn={0012-7094},
   doi={10.1215/00127094-2009-062},
}

		
\bib{FOOO:II}{article}{
   author={Fukaya, Kenji},
   author={Oh, Yong-Geun},
   author={Ohta, Hiroshi},
   author={Ono, Kaoru},
   title={Lagrangian Floer theory on compact toric manifolds II: bulk
   deformations},
   journal={Selecta Math. (N.S.)},
   volume={17},
   date={2011},
   number={3},
   pages={609--711},
   issn={1022-1824},
   doi={10.1007/s00029-011-0057-z},
}

\bib{FOOO:Ast}{article}{
   author={Fukaya, Kenji},
   author={Oh, Yong-Geun},
   author={Ohta, Hiroshi},
   author={Ono, Kaoru},
   title={Lagrangian Floer theory and mirror symmetry on compact toric
   manifolds},
   language={English, with English and French summaries},
   journal={Ast\'{e}risque},
   number={376},
   date={2016},
   pages={vi+340},
   issn={0303-1179},
   isbn={978-2-85629-825-1},
}

\bib{GGI}{article}{
   author={Galkin, Sergey},
   author={Golyshev, Vasily},
   author={Iritani, Hiroshi},
   title={Gamma classes and quantum cohomology of Fano manifolds: gamma
   conjectures},
   journal={Duke Math. J.},
   volume={165},
   date={2016},
   number={11},
   pages={2005--2077},
   issn={0012-7094},
   doi={10.1215/00127094-3476593},
}

\bib{Galkin-Iritani}{article}{
author={Galkin, Sergey}, 
author={Iritani, Hiroshi}, 
title={Gamma conjecture via mirror symmetry},
note={\arxiv{1508.00719}} 
}

\bib{GKZ}{article}{
   author={Gel\cprime fand, I. M.},
   author={Zelevinski\u\i , A. V.},
   author={Kapranov, M. M.},
   title={Hypergeometric functions and toric varieties},
   language={Russian},
   journal={Funktsional. Anal. i Prilozhen.},
   volume={23},
   date={1989},
   number={2},
   pages={12--26},
   issn={0374-1990},
   translation={
      journal={Funct. Anal. Appl.},
      volume={23},
      date={1989},
      number={2},
      pages={94--106},
      issn={0016-2663},
   },
   doi={10.1007/BF01078777},
}

\bib{Givental:ICM}{article}{
   author={Givental, Alexander B.},
   title={Homological geometry and mirror symmetry},
   conference={
      title={ 2},
      address={Z\"urich},
      date={1994},
   },
   book={
      publisher={Birkh\"auser},
      place={Basel},
   },
   date={1995},
   pages={472--480},
}


\bib{Givental:mirror_thm_toric}{article}{
   author={Givental, Alexander},
   title={A mirror theorem for toric complete intersections},
   conference={
      title={Topological field theory, primitive forms and related topics},
      address={Kyoto},
      date={1996},
   },
   book={
      series={Progr. Math.},
      volume={160},
      publisher={Birkh\"auser Boston, Boston, MA},
   },
   date={1998},
   pages={141--175},
}

\bib{Golyshev:classification}{article}{
   author={Golyshev, Vasily V.},
   title={Classification problems and mirror duality},
   conference={
      title={Surveys in geometry and number theory: reports on contemporary
      Russian mathematics},
   },
   book={
      series={London Math. Soc. Lecture Note Ser.},
      volume={338},
      publisher={Cambridge Univ. Press, Cambridge},
   },
   date={2007},
   pages={88--121},
   doi={10.1017/CBO9780511721472.004},
}

\bib{Golyshev:techniques}{article}{
author={Golyshev, Vasily}, 
title={Techniques to compute monodromy of differential 
equations of mirror symmetry},  
note={to appear in the same volume.} 
}

\bib{Golyshev-Zagier}{article}{
   author={Golyshev, V. V.},
   author={Zagier, D.},
   title={Proof of the gamma conjecture for Fano 3-folds with a Picard
   lattice of rank one},
   language={Russian, with Russian summary},
   journal={Izv. Ross. Akad. Nauk Ser. Mat.},
   volume={80},
   date={2016},
   number={1},
   pages={27--54},
   issn={1607-0046},
   translation={
      journal={Izv. Math.},
      volume={80},
      date={2016},
      number={1},
      pages={24--49},
      issn={1064-5632},
   },
   doi={10.4213/im8343},
}

\bib{Gonzalez-Woodward:tmmp}{article}{
author={Gonz\'alez, Eduardo}, 
author={Woodward, Chris}, 
title={Quantum cohomology and toric minimal model programs},  
year={2012}, 
note={\arxiv{1010.2118}}
}

\bib{Hori-Romo}{article}{
author={Hori, Kentaro},
author={Romo, Mauricio}, 
title={Exact results in two-dimensional (2,2) supersymmetric gauge theory
with boundary}, 
note={\arxiv{1308.2438}} 
}


\bib{Hori-Vafa}{article}{
author={Hori, Kentaro}, 
author={Vafa, Cumrum}, 
title={Mirror symmetry}, 
note={\arxiv{hep-th/0002222}}, 
year={2000}
}

\bib{Hosono:centralcharges}{article}{
   author={Hosono, Shinobu},
   title={Central charges, symplectic forms, and hypergeometric series in
   local mirror symmetry},
   conference={
      title={Mirror symmetry. V},
   },
   book={
      series={AMS/IP Stud. Adv. Math.},
      volume={38},
      publisher={Amer. Math. Soc., Providence, RI},
   },
   date={2006},
   pages={405--439},
}

\bib{HKTY:CICY}{article}{
   author={Hosono, S.},
   author={Klemm, A.},
   author={Theisen, S.},
   author={Yau, S.-T.},
   title={Mirror symmetry, mirror map and applications to complete
   intersection Calabi-Yau spaces},
   journal={Nuclear Phys. B},
   volume={433},
   date={1995},
   number={3},
   pages={501--552},
   issn={0550-3213},
   doi={10.1016/0550-3213(94)00440-P},
}

\bib{Iritani:Integral}{article}{
   author={Iritani, Hiroshi},
   title={An integral structure in quantum cohomology and mirror symmetry
   for toric orbifolds},
   journal={Adv. Math.},
   volume={222},
   date={2009},
   number={3},
   pages={1016--1079},
   issn={0001-8708},
   doi={10.1016/j.aim.2009.05.016},
}

\bib{Iritani:big_equivariant}{article}{
author={Iritani, Hiroshi}, 
title={A mirror construction for the big equivariant quantum cohomology 
of toric manifolds}, 
journal={Math. Ann.}, 
year={2017}, 
volume={368}, 
pages={279--316}
}



\bib{KKP:Hodge}{article}{
   author={Katzarkov, L.},
   author={Kontsevich, M.},
   author={Pantev, T.},
   title={Hodge theoretic aspects of mirror symmetry},
   conference={
      title={From Hodge theory to integrability and TQFT tt*-geometry},
   },
   book={
      series={Proc. Sympos. Pure Math.},
      volume={78},
      publisher={Amer. Math. Soc., Providence, RI},
   },
   date={2008},
   pages={87--174},
   doi={10.1090/pspum/078/2483750},
}

\bib{KKP:Bogomolov}{article}{
   author={Katzarkov, Ludmil},
   author={Kontsevich, Maxim},
   author={Pantev, Tony},
   title={Bogomolov-Tian-Todorov theorems for Landau-Ginzburg models},
   journal={J. Differential Geom.},
   volume={105},
   date={2017},
   number={1},
   pages={55--117},
   issn={0022-040X},
}

\bib{LLY:I}{article}{
   author={Lian, Bong H.},
   author={Liu, Kefeng},
   author={Yau, Shing-Tung},
   title={Mirror principle. I},
   journal={Asian J. Math.},
   volume={1},
   date={1997},
   number={4},
   pages={729--763},
   issn={1093-6106},
   doi={10.4310/AJM.1997.v1.n4.a5},
}

\bib{LLY:II}{article}{
   author={Lian, Bong H.},
   author={Liu, Kefeng},
   author={Yau, Shing-Tung},
   title={Mirror principle. II},
   note={Sir Michael Atiyah: a great mathematician of the twentieth
   century},
   journal={Asian J. Math.},
   volume={3},
   date={1999},
   number={1},
   pages={109--146},
   issn={1093-6106},
   doi={10.4310/AJM.1999.v3.n1.a6},
}

\bib{Libgober}{article}{
   author={Libgober, Anatoly},
   title={Chern classes and the periods of mirrors},
   journal={Math. Res. Lett.},
   volume={6},
   date={1999},
   number={2},
   pages={141--149},
   issn={1073-2780},
   doi={10.4310/MRL.1999.v6.n2.a2},
}

\bib{Mann-deGregorio}{article}{
   author={de Gregorio, Ignacio},
   author={Mann, \'{E}tienne},
   title={Mirror fibrations and root stacks of weighted projective spaces},
   journal={Manuscripta Math.},
   volume={127},
   date={2008},
   number={1},
   pages={69--80},
   issn={0025-2611},
   doi={10.1007/s00229-008-0185-8},
}

\bib{Mann-Reichelt}{article}{
author ={Mann, Etienne},  
author ={Reichelt, Thomas}, 
title ={Logarithmic degeneration of Landau-Ginzburg models 
for toric orbifolds and global $tt^*$-geometry}, 
note ={\arxiv{1605.08937}} 
}

\bib{McDuff-Tolman}{article}{
   author={McDuff, Dusa},
   author={Tolman, Susan},
   title={Topological properties of Hamiltonian circle actions},
   journal={IMRP Int. Math. Res. Pap.},
   date={2006},
   pages={72826, 1--77},
   issn={1687-3017},
}

\bib{Mochizuki:twistor_GKZ}{article}{
author={Mochizuki, Takuro}, 
title={Twistor property of GKZ-hypergeometric systems},
note={\arxiv{1501.04146}} 
}

\bib{Reichelt-Sevenheck}{article}{
   author={Reichelt, Thomas},
   author={Sevenheck, Christian},
   title={Logarithmic Frobenius manifolds, hypergeometric systems and
   quantum $\mathcal{D}$-modules},
   journal={J. Algebraic Geom.},
   volume={24},
   date={2015},
   number={2},
   pages={201--281},
   issn={1056-3911},
   doi={10.1090/S1056-3911-2014-00625-1},
}

\bib{Rietsch:GmodP}{article}{
   author={Rietsch, Konstanze},
   title={A mirror symmetric construction of $qH^\ast_T(G/P)_{(q)}$},
   journal={Adv. Math.},
   volume={217},
   date={2008},
   number={6},
   pages={2401--2442},
   issn={0001-8708},
   doi={10.1016/j.aim.2007.08.010},
}

\bib{Ritter}{article}{
   author={Ritter, Alexander F.},
   title={Circle actions, quantum cohomology, and the Fukaya category of
   Fano toric varieties},
   journal={Geom. Topol.},
   volume={20},
   date={2016},
   number={4},
   pages={1941--2052},
   issn={1465-3060},
   doi={10.2140/gt.2016.20.1941},
}


\bib{Saito:higher_residue}{article}{
   author={Saito, Kyoji},
   title={The higher residue pairings $K_{F}^{(k)}$\ for a family of
   hypersurface singular points},
   conference={
      title={Singularities, Part 2},
      address={Arcata, Calif.},
      date={1981},
   },
   book={
      series={Proc. Sympos. Pure Math.},
      volume={40},
      publisher={Amer. Math. Soc., Providence, RI},
   },
   date={1983},
   pages={441--463},
}

\bib{Seidel}{article}{
author={Seidel, Paul}, 
title={$\pi_1$ of symplectic automorphism groups and invertibles 
in quantum homology rings}, 
journal={Geom. Funct. Anal.}, 
volume={7}, 
number={6}, 
pages={1046--1095}, 
year={1997} 
}

\bib{Smith}{article}{
author={Smith, Jack}, 
title={Quantum cohomology and closed-string mirror symmetry 
for toric varieties},
note={\arxiv{1802.00424}},
year={2018}
}

\bib{SYZ}{article}{
   author={Strominger, Andrew},
   author={Yau, Shing-Tung},
   author={Zaslow, Eric},
   title={Mirror symmetry is $T$-duality},
   journal={Nuclear Phys. B},
   volume={479},
   date={1996},
   number={1-2},
   pages={243--259},
   issn={0550-3213},
   doi={10.1016/0550-3213(96)00434-8},
}


\end{biblist}
\end{bibdiv}
\providecommand{\arxiv}[1]{\href{http://arxiv.org/abs/#1}{arXiv:#1}}

\end{document}